\documentclass[11pt]{article}

\usepackage{amssymb}
\usepackage{graphics}
\usepackage{amsmath,enumerate}
\usepackage{enumitem}
\newtheorem{thm}{Theorem}
\newtheorem{clm}{Claim}

\def \no {\noindent}
\def \es {\emptyset}
 \def \sm {\setminus}

\title{Coloring ($P_6$, diamond, $K_4$)-free graphs}

\author{T. Karthick\thanks{Corresponding author. Computer Science Unit, Indian Statistical
Institute, Chennai Centre, Chennai-600113, India.
E-mail: {karthick@isichennai.res.in}} \and Suchismita Mishra\thanks{Department of Mathematics, Indian Institute of Technology Madras, Chennai-600036, India.}}

\begin{document}
\maketitle

\begin{abstract}
We show that every ($P_6$, diamond, $K_4$)-free graph  is $6$-colorable. Moreover, we give an example of a ($P_6$, diamond, $K_4$)-free graph $G$ with $\chi(G) = 6$. This generalizes some known results in the literature.
\end{abstract}

\section{Introduction}
We consider simple, finite, and undirected graphs. For notation and terminology  not
defined here we refer to \cite{West}.  Let $P_n$,
$C_n$, $K_n$ denote the induced path, induced cycle and complete graph
on $n$ vertices respectively.  If $G_{1}$ and $G_{2}$ are two vertex disjoint graphs, then
  their {\it union} $G_1\cup G_2$ is the graph with $V(G_1\cup G_2)$
  $=V(G_1)\cup V(G_2)$ and $E(G_{1}\cup G_{2})$ $=$ $E(G_1)\cup
  E(G_2)$. Similarly, their {\it join} $G_{1}+G_{2}$ is the graph with
  $V(G_1+G_2)$ $=$ $V(G_1)\cup V(G_2)$ and $E(G_1 + G_2)$
  $=$ $E(G_1)\cup E(G_2)$$\cup \{(x,y)\mid x\in V(G_1), ~y\in
  V(G_2)\}$.    For any positive integer $k$, $kG$ denotes the union
  of $k$ graphs each isomorphic to $G$.  If $\cal{F}$ is a family of
graphs, a graph $G$ is said to be \emph{$\cal{F}$-free} if it contains no
induced subgraph isomorphic to any member of $\cal{F}$.  A \emph{clique} (independent set) in a
graph $G$ is a set of vertices that are pairwise adjacent (non-adjacent) in $G$. The \emph{clique number} of $G$, denoted by
$\omega(G)$, is the size of a maximum clique in $G$.

 A  \emph{$k$-coloring} of a graph $G = (V, E)$ is a mapping $f : V \rightarrow \{1,2, \ldots, k\}$
such that $f(u) \neq f(v)$ whenever $uv \in E$. We say that $G$ is \emph{$k$-colorable} if $G$ admits a $k$-coloring. That is, a partition of the vertex set $V(G)$ into $k$ independent sets. The
\emph{chromatic number} of $G$, denoted by $\chi(G)$, is the smallest positive integer $k$ such
that $G$ is $k$-colorable. It is well known that a graph  is $2$-colorable if and only if it is bipartite.  Given an integer $k$, the $k$-{\sc Coloring} problem is that of testing whether a given graph is $k$-colorable.   The $k$-{\sc Coloring} problem is
$NP$-complete for every fixed $k\ge 3$ \cite{GJS, Karp}.
   The problem of finding the maximum chromatic number of graphs without forbidden induced
subgraphs from some finite/infinite set and with small clique number  is well studied and still receives much attention. We refer to \cite{Randerath-Schiermeyer-survey} for a survey and we give some of them that are not in \cite{Randerath-Schiermeyer-survey}.  We note that some of the cited results are consequences of much stronger results available in the literature.  Mycielski \cite{Mycielski} showed that for any  integer $k$, there exists a triangle-free graph with chromatic number $k$.  Fan et al. \cite{FXYY} showed that every (fork, $K_3$)-free graph with odd-girth at least $7$ is $3$-colorable.  Pyatkin \cite{Pyatkin} showed that every ($2P_3, K_3$)-free graph is 4-colorable.
 Esperet et al. \cite{ELMM} showed that every ($P_5, K_4$)-free graph is $5$-colorable.
     It follows from a result of Gravier, Ho\'ang and Maffray \cite{GHM} that  every ($P_6, K_3$)-free graph is $4$-colorable (see also \cite{Randerath-Schiermeyer-survey}), and that every ($P_6, K_4$)-free graph is $16$-colorable. Randerath et al.  \cite{BST} showed that every ($P_6, K_{1,3}, W_5,  DD, K_4$)-free graph is $3$-colorable, where $DD$ denotes the double-diamond graph, and $W_5$ is the 5-wheel.  It follows from a result of \cite{CK} that every ($P_2\cup P_3, C_4, K_4$)-free graph is $4$-colorable, and from a result of \cite{KM} that every ($P_5$, diamond, $K_4$)-free graph is $4$-colorable. Chudnovsky et al. \cite{CSRT-K4-free} showed that every   (odd hole, $K_4$)-free graph is 4-colorable. This implies that every ($P_6, C_5, K_4$)-free graph is $4$-colorable.  Addario-Berry et al. \cite{ACHRS}  showed that every  (even hole, $K_4$)-free graph is 5-colorable, and Kloks, M\"uller and Vuskovic \cite{KMV} showed that every (even-hole, diamond, $K_4$)-free is $4$-colorable.

\begin{figure}
\centering
 \includegraphics{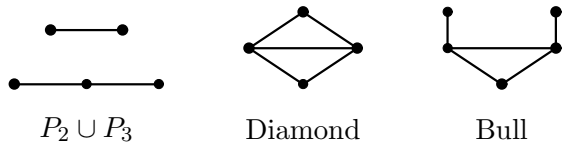}
\caption{Some special graphs.}\label{special}
\end{figure}

In this note, we first give an alternative and simple proof  to the fact that every ($P_2 \cup P_3$,  diamond, $K_4$)-free graph is $6$-colorable given in \cite{AC}. Then we show that the conclusion holds even for a general class of graphs, namely ($P_6$,  diamond, $K_4$)-free graphs. That is, we show that every ($P_6$,  diamond, $K_4$)-free graph is $6$-colorable. Moreover, we give an example of a ($P_6$,  diamond, $K_4$)-free graph $G$ with $\chi(G) =6$. This generalizes the aforementioned results for ($P_5$, diamond, $K_4$)-free graphs, ($P_6, K_3$)-free graphs and  ($P_2 \cup P_3$,  diamond, $K_4$)-free graphs. The proof of our results depend on a sequence of partial results  given below, and  we give tight examples for each of them. See Figure~\ref{special}.
\begin{itemize}
\item Let $G$ be a ($P_2 \cup P_3$,  diamond, $K_4$)-free graph that contains a non-dominating $K_3$.
  Then $G$ is $4$-colorable.
\item  Let $G$ be a ($P_2 \cup P_3$,   diamond, $K_4$)-free graph such that every triangle in $G$ dominates $G$.
  Then $G$ is $6$-colorable.
\item  Let $G$ be a ($P_6$,  diamond, bull, $K_4$)-free graph. Then $G$ is 4-colorable.
\item  Let $G$ be a ($P_6$, diamond, $K_4$)-free graph that contains an induced bull. Then $G$ is $6$-colorable.
    \end{itemize}

Note that the class of ($H$, diamond)-free graphs, for various $H$, is well studied in variety of contexts in the literature.  Chudnovsky et al. \cite{CGSZ} showed that there are exactly six 4-critical
($P_6$, diamond)-free graphs.
Tucker \cite{Tucker} gave an $O(kn^2)$ time algorithm for $k$-{\sc Coloring}  perfect diamond-free graphs. It is also known that $k$-{\sc Coloring} is polynomial-time solvable for (even-hole, diamond)-free graphs \cite{KMV}  as well as for  (hole, diamond)-free graphs \cite{BGM}. Dabrowski et al. \cite{DDP} showed that if $H$ is a graph on at most five vertices, then $k$-{\sc Coloring} is polynomial time
solvable for ($H$, diamond)-free graphs, whenever $H$ is a linear forest and NP-complete otherwise. However, the computational complexity of the
$k$-{\sc Coloring} problem for ($P_6$, diamond)-free graphs is open. It is also known that the {\sc Maximum Weight Independent Set} problem is solvable in polynomial time for ($P_6$, diamond)-free graphs \cite{Mosca} as well as for (hole, diamond)-free graphs \cite{BGM}.

We devote the rest of this section to the notations and terminologies used in this paper.
For any integer $k$, we simply write $[k]$ to denote the set $\{1, 2, \ldots, k\}$. Let $G$ be a graph, with vertex-set $V(G)$ and edge-set $E(G)$. A \emph{diamond} or a  $K_4-e$ is the graph with vertex set $\{a, b, c, d\}$ and edge set $\{ab, bc, cd, ad, bd\}$.

A graph $G$ is called {\it perfect} if $\chi(H)= \omega(H)$, for every
 induced subgraph $H$ of $G$.   By the \emph{Strong Perfect Graph Theorem}
\cite{Chudnovsky}, A graph is perfect if and only if it
contains no odd hole (chordless cycle) of length at least $5$ and no
odd antihole (complement graph of a hole) of length at least $5$.

For $x \in V(G)$, $N(x)$ denotes the set of all neighbors of $x$ in
$G$. The neighborhood
$N(X)$ of a subset $X \subseteq V(G)$ is the set $\{u \in
V(G)\setminus X \mid u \mbox{~is adjacent to a vertex of }X\}$. For any two disjoint subsets $S,~T \subseteq V(G)$,
  $[S, T]$ denotes the edge-set  $\{uv \in E(G)\mid u \in S,\ v \in T\}$. The set $[S, T]$ is said to be {\it complete} if
every vertex in $S$ is adjacent to every vertex in $T$. Also, for $S \subseteq V(G)$, let $G[S]$ denotes the subgraph induced by
  $S$ in $G$, and for convenience we simply write $[S]$ instead of $G[S]$. The length of a path is the number of edges in it. The length of a
shortest path between two vertices $x$ and $y$ is denoted by $dist(x,y)$. For $S \subseteq V(G)$ and $x \in V(G) \setminus S$, we define dist$(x, S):= \min\{\mbox{dist}(x, y) \mid y \in S\}$.
  We say that a subgraph $H$ of $G$ is \emph{dominating}  if every vertex in $V(G)\sm V(H)$ has a neighbor in $H$; otherwise, it is a \emph{non-dominating} subgraph. We say that a graph $H$ is obtained from
$G$ by \emph{duplication} if it can be obtained from $G$ by substituting
independent sets for some of the vertices in $G$.

\section{Our Results}

We use the following preliminary results often. Let $G$ be a (diamond, $K_4$)-free graph. Then  the following hold:

\begin{enumerate}
\item[(R1)] If $T$ is a triangle in $G$, then any vertex $p \in V(G) \sm V(T)$ has at most one neighbor in $T$.
\item[(R2)] For any $v \in V(G)$, $N(v)$ induces a $P_3$-free graph, and hence $[N(v)]$ is a union of $K_2$'s and $K_1$'s.
 \item[(R3)] For any two non-adjacent vertices $x$ and $y$ in $G$, the set of common neighbors of $x$ and $y$ is an independent set.
  \item[(R4)] For any two adjacent vertices $x$ and $y$ in $G$, the number of common neighbors of $x$ and $y$ is at most one.

\item[(R5)] Let $G$ be a diamond-free graph with $n$ vertices. Let $v_1, v_2, \ldots, v_k$ $(1 \leq k \leq n)$ be vertices in $G$ such
that  $N(v_i)$ is an independent set, for each $i \in \{1, 2 \ldots, k\}$. Then the duplicated graph $H$ obtained from $G$ by substituting an independent set for each $v_i$, $i \in \{1, 2 \ldots, k\}$ is also diamond-free.
\end{enumerate}

\begin{thm}\label{non-dom-k3-4-col}
  Let $G$ be a connected ($P_2 \cup P_3$,  diamond, $K_4$)-free graph that contains a non-dominating $K_3$.
  Then $G$ is $4$-colorable.
\end{thm}
\no{\it Proof of Theorem~$\ref{non-dom-k3-4-col}$}. Let $T$ be a non-dominating triangle in $G$ induced by the vertex-set
$N_0:= \{v_1, v_2, v_3\}$.
 For $j \geq 1$,  let $N_j$ denote the set $\{y \in V(G)\setminus N_0 \mid d(y, N_0) = j\}$.
  Then by (R1), every vertex of $N_1$
has at most one neighbor in $N_0$; more precisely:

\begin{enumerate}
\item[(1)] If $x \in N_1$, then $[N(x)\cap N_0]$ is isomorphic to $K_1$.
\end{enumerate}
For $i \in [3]$, let $A_i:= \{x\in N_1  \mid  N(x)\cap N_0 =\{v_i\}\}$. Then since $G$ is $K_4$-free,  by (R2), each $[A_i]$ is a union of $K_2$'s and $K_1$'s, and so $[A_i]$ is bipartite. Let $(A_1', A_1'')$ be a bipartition of $A_1$ such that $A_1'$ is a maximal independent set of $A_1$.

\begin{enumerate}
\item[(2)]$N_j = \emptyset$, for all $j \geq 3$.

{\it Proof}. Suppose not, and let $y\in N_3$. Then by the definition of $N_3$, there exist vertices $x_1 \in N_1$ and $x_2 \in N(x)\cap N_2$ such that $x_1-x_2-y$ is a $P_3$ in $G$. By (1), $x_1 \in A_i$, for some $i \in [3]$, say $i =1$. But then $\{v_2, v_3, x_1, x_2, y\}$ induces a $P_2\cup P_3$ in $G$, which is a contradiction. $\diamond$
\end{enumerate}

Now, since $T$ is non-dominating, $N_2 \neq \emptyset$. Also, since $G$ is ($P_2 \cup P_3$)-free,
$[N_2]$ is a $P_3$-free graph, and hence  $[N_2]$ is a union of complete graphs. Moreover:
\begin{enumerate}
\item[(3)] If $C$ is a component of $[N_2]$, then there exists $x \in N_1$ such that $[N(x) \cap N_2] \cong C$.

{\it Proof}. Otherwise, $G$ induces a $P_2 \cup P_3$ or diamond, a contradiction. $\diamond$
\end{enumerate}

 Furthermore, since $G$ is $K_4$-free, by (3), $[N_2]$ is a union of $K_2$'s and $K_1$'s.
In particular, $[N_2]$ is bipartite. Let $(N_2', N_2'')$ be a bipartition of $N_2$.

\begin{enumerate}
\item[(4)]
  Let $i, j \in [3]$. Let $x\in A_i$ be  such that $N(x)\cap N_2 \neq \es$. Then $A_j$ ($j \neq i$) is  an independent set.

\no{\it Proof}. We may assume that $i =1$. Let $z \in N(x) \cap N_2$. Suppose to the contrary that there are vertices $y_1, y_2 \in A_2$ such that $y_1y_2 \in E$. Then since $\{y_1, y_2, v_2, x\}$ does not induce a diamond in $G$, we have either $xy_1 \notin E$ or $xy_2 \notin E$. Assume $xy_1 \notin E$. Now, since $\{x, z, y_1, v_2, v_3\}$ does not induce a $P_2 \cup P_3$ in $G$, $zy_1 \in E$. Then since $\{v_1, v_3, z, y_1, y_2\}$ does not induce a $P_2 \cup P_3$ in $G$, $zy_2 \in E$. But, then $\{v_2, y_1, y_2, z\}$ induces a diamond in $G$, a contradiction. So, $A_2$ is an independent set. Similarly, $A_3$ is also an independent set.  $\diamond$
\end{enumerate}

Now,  we show that $G$ is 4-colorable using the above properties.
 Suppose that $K_2 \sqsubseteq [N_2]$. Then by (3), there is a vertex $x \in A_i$ such that $[N(x)\cap N_2] \cong K_2$.
 We may assume that $i =1$, and let $y_1, y_2 \in N(x) \cap N_2$ be such that $y_1y_2 \in E$. Then $A_2 = \emptyset = A_3$. Otherwise, if $w \in A_2\cup A_3$, then since $\{w, v_2, v_3, y_1, y_2\}$ or $\{v_1, v_3, w, y_1, y_2\}$ or $\{v_1, v_2, w, y_1, y_2\}$ does not induce a $P_2 \cup P_3$ in $G$, we have  $wy_1, wy_2 \in E$. But, then $\{x, w, y_1, y_2\}$ induces a diamond or a $K_4$ in $G$, a contradiction.  Then we define the sets $S_1:= \{v_1\}\cup N_2', S_2:= \{v_2\}\cup A_1', S_3 := \{v_3\}\cup A_1''$, and $S_4:= N_2''$. Clearly, $S_j$ is an independent set, for each $j \in \{1, 2, 3, 4\}$. Hence $(S_1, S_2, S_3, S_4)$ is a $4$-coloring of $G$.

So, assume that $N_2$ $(\neq \es)$ is an independent set. Then there exists $x \in A_i$ such that $N(x)\cap N_2 \neq \emptyset$.  We may assume that $i=1$, and let $z \in N(x) \cap N_2$. By (4), $A_2$ and $A_3$ are independent sets.
Moreover, if $y\in A_1''$, then $N(y) \cap N_2 = \es$. (Otherwise, since $A_1'$ is maximal, there exist $y' \in A_1'$ and $z' \in N(y) \cap N_2$ such that
$yy', yz' \in E$. Then since $\{v_2, v_3, y, y', z'\}$ does not induces a $P_2\cup P_3$ in $G$, $y'z' \in E$. But, then $\{v_1, y, y', z'\}$ induces a diamond in $G$, a contradiction.) Now, we define $S_1:=\{v_1\}\cup A_2, S_2:= \{v_2\}\cup A_1', S_3 := \{v_3\}\cup A_1''\cup N_2,$ and $S_4:= A_3$. Then $S_j$ is an independent set, for each $j \in \{1, 2, 3, 4\}$. Hence $(S_1, S_2, S_3, S_4)$ is a $4$-coloring of $G$.

Thus, Theorem \ref{non-dom-k3-4-col} is proved. \hfill{$\Box$}

\begin{figure}
\centering
 \includegraphics{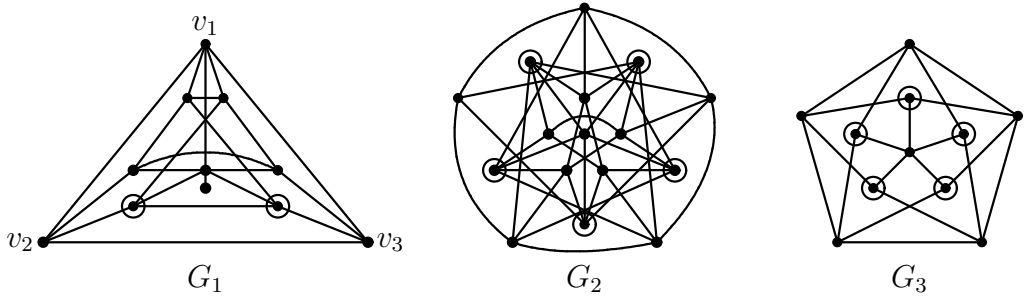}
\caption{Some extremal graphs.}\label{ext}
\end{figure}

\medskip
The bound given in Theorem~\ref{non-dom-k3-4-col} is tight. For example, consider the duplicated graph $G$ obtained from $G_1$ (shown in Figure~\ref{ext}) by substituting each vertex indicated in circle by an independent set (of order $\geq 1$).  As the graph $G$ is highly symmetric, using (R5), there are not too many cases to directly verify that $G$ is ($P_2 \cup P_3$,  diamond, $K_4$)-free. Also, it is easy to see that $G$ contains a non-dominating triangle $T$ with vertices $\{v_1, v_2, v_3\}$, and  that $\chi(G) =4$.

\begin{thm}\label{dom-k3-6-col}
  Let $G$ be a connected ($P_2 \cup P_3$,  diamond,  $K_4$)-free graph such that every triangle in $G$ dominates $G$.
  Then $G$ is $6$-colorable.
\end{thm}
\no{\it Proof of Theorem~$\ref{dom-k3-6-col}$}.  Let $K$ be a $K_3$ in $G$ with vertices $\{x_1, x_2, x_3\}$ that dominates $G$. Since $G$ is ($K_4$, diamond)-free and $K$ dominates $G$, every vertex in  $V(G) \setminus V(K)$
has exactly one neighbor in $K$. For $i \in [4]$, let $A_i:= \{x\in V(G) \setminus V(K)  \mid  N(x)\cap V(K) =\{x_i\}\}$.  Then by (R2), each $[A_i]$ is a union of $K_2$'s and $K_1$'s. Now, it is easy to check that $G$ is $6$-colorable.   \hfill{$\Box$}

\medskip
The bound given in Theorem~\ref{dom-k3-6-col} is tight. For example, consider the graph $G$ which is isomorphic to the complement of the 16-regular \emph{Schl\"afli graph} on 27 vertices. It is verified that $G$ is ($P_2 \cup P_3$, diamond)-free and it is well known that $\chi(G) = 6$ and $\omega(G) = 3$; see (https://hog.grinvin.org/ViewGraphInfo.action?id=19273).

\begin{thm}
Let $G$ be a ($P_2 \cup P_3$,  diamond, $K_4$)-free graph. Then $G$ is $6$-colorable.
\end{thm}
\no{\it Proof}. Follows by Theorems~\ref{non-dom-k3-4-col} and \ref{dom-k3-6-col}. \hfill{$\Box$}

\begin{thm} \label{PDKB-4-col}
 Let $G$ be a ($P_6$, diamond, bull, $K_4$)-free graph. Then $G$ is 4-colorable.
\end{thm}

\no{\it Proof}. If $G$ is perfect, then $G$ is $3$-colorable and the theorem holds. So we may assume that $G$ is connected and $G$ is not perfect. Since $G$ is $P_6$-free, $G$ contains no hole of length at least $7$, and since $G$ is diamond-free, $G$ contain no anti-hole of length at least $7$.
Thus, it follows from the Strong Perfect Graph Theorem \cite{Chudnovsky} that $G$ contains a $5$-hole (hole of length $5$), say $N_0$ with vertex-set
$\{v_1, v_2, v_3, v_4, v_5\}$, and edge-set $\{v_1v_2, v_2v_3, v_3v_4, v_4v_5, v_5v_1\}$. Throughout this proof, we take all
the subscripts of $v_i$ to be modulo $5$. For any integer $j \geq 1$, let $N_j$ denote the set $\{x \in V(G) \mid d(x, N_0) = j \}$.

\begin{enumerate}
\item[(1)] If $x\in N_1$, then $[N(x)\cap N_0]$ is isomorphic to either $K_1$ or $2K_1$.

\no{\it Proof}. Suppose not. Then there exists an $i\in [5]$ such that $\{v_i, v_{i+1}\} \subseteq N(x)$. Then since $\{v_{i-1}, v_i,  v_{i+1}, v_{i+2}, x\}$ does not induce a bull in $G$, either $v_{i-1} \in N(x)$ or $v_{i+2} \in N(x)$. But then either $\{v_i, v_{i+1},  v_{i+2}, x\}$ or $\{v_{i-1}, v_{i}, v_{i+1}, x\}$  induces a diamond in $G$, a contradiction. So $(1)$ holds. $\diamond$
\end{enumerate}

\no{}By (1), we partition $N_1$ as follows:
For any $i\in [5]$, $i$ mod $5$, let:
\begin{eqnarray*}
W_{i} &:=& \{ x\in N_1 \mid N(x)\cap N_0 = \{ v_i \} \},\\
Y_i &:= & \{x \in N_1 \mid N(x)\cap N_0 = \{v_{i-1},v_{i+1}\}.
\end{eqnarray*}
Moreover, let $W := W_1 \cup \cdots \cup W_5$ and $Y := Y_1 \cup \cdots \cup Y_5$.

\begin{enumerate}
\item[(2)] For   $i \in [5]$, $i$ mod $5$, the following hold:
\begin{enumerate}[label=(\roman*)]
\item $[W_i]$ is union of $K_2$'s and $K_1$'s.

\item $[W_i, W_{i+1}]=\es$.
%\item[(iii)] $[W_i, W_{i+2}]$ and $[W_i, W_{i-2}]$ are complete.

\item $[W_i, W_{i+2}]$ is complete. In particular, if $W_i\neq \es$, then $W_{i+2}$ and $W_{i-2}$ are independent sets.
\item $Y_i \cup Y_{i+2}$ is an independent set.
\item $[W_i, Y_{i+1}] = \es$
\end{enumerate}
\end{enumerate}
\no{\it Proof of $(2)$}. We prove for $i = 1$.

\smallskip
\no{$(i)$}:  Follows by the definition of $W_1$, and by (R2).

\smallskip
\no{$(ii)$}: Suppose not. Then there exist vertices $x\in W_1$ and $y \in W_{2}$ such that $xy \in E$. But, then
$\{x, y, v_{2}, v_{3}, v_{4}, v_{5}\}$  induces a $P_6$ in $G$, which is a contradiction. So (ii) holds.

\smallskip
\no{$(iii)$}: Let $x \in W_1$. We show that $W_3$ is an independent set.  If not,  then there exist adjacent vertices $y$ and $z$ in $W_3$. Since $\{x, v_1, v_{5}, v_{4}, v_{3}, y\}$ and $\{x, v_1, v_{5}, v_{4}, v_{3}, z\}$ do not induce a $P_6$ in $G$, we have $xy \in E$ and $xz \in E$. But, then
$\{x, y, z, v_3\}$ induces a diamond in $G$, which is a contradiction. So, $W_3$ is an independent set. Similarly,  $W_4$ is also an independent set.
Thus (iii) holds.

\smallskip
\no{$(iv)$}: Suppose not. Then there exist adjacent vertices $x$ and $y$ in $Y_1 \cup Y_{3}$. If both $x$ and $y$ are in $Y_1$ or in $Y_3$, then
either $\{v_1, x, v_5, y\}$ or $\{v_2, x, v_4, y\}$ induces a diamond in $G$, a contradiction. So, we may assume that $x \in Y_1$ and $y \in Y_3$.
But, then $\{v_1, v_2, y, v_4, x\}$ induces a bull in $G$, which is a contradiction. So (iv) holds.

\smallskip
\no{$(v)$}:   Suppose not. Then there exist vertices $x\in W_1$ and $y \in Y_{2}$ such that $xy \in E$. But, then $\{v_5, v_1, y, v_3, x\}$ induces a bull in $G$, which is a contradiction. So (v) holds. $\diamond$

\medskip
By (2)(i), for each $i\in [5]$, $W_i$ is bipartite. Let $(W_i',W_i'')$ is a bipartition of $W_i$.

\begin{enumerate}
\item[(3)] The following hold:
\begin{enumerate}[label=(\roman*)]
\item  $[W, N_2] = \es$.
\item If $xx'$ is an edge in $[N_2]$, then $N(x)\cap N_1=N(x')\cap N_1$. Moreover, $|N(x)\cap N_1|=1$.
\item If $C$ is a component of $[V(G)\sm (N_0\cup N_1)]$, then there exists a vertex $y \in Y$ such that $N(y)\cap V(C) = V(C)$.

\end{enumerate}
\end{enumerate}

\no{\it Proof of $(3)$}. $(i)$:  Suppose not. Then there exist vertices $x\in W$ and $y \in N_{2}$ such that $xy \in E$. We may assume that $x \in W_1$. But, then  $\{y, x, v_1, v_{2},v_{3}, v_{4}\}$ induces a $P_6$ in $G$ which is a contradiction. So (i) holds.

\smallskip
\no{$(ii)$}:  Let $y \in N(x)\cap N_1$. Then by (i),  $y \in Y_i$, for some $i$. Say $y \in Y_1$. Since $\{x', x, y, v_{2},v_{3}, v_{4}\}$ does not induce a $P_6$ in $G$, $x'y\in E$. So, $N(x)\cap N_1 \subseteq N(x')\cap N_1$. Similarly, $N(x')\cap N_1\subseteq N(x)\cap N_1$. Hence, $N(x)\cap N_1 = N(x')\cap N_1$. Moreover,  if $|N(x)\cap N_1| > 1$, and if  $y_1, y_2 \in N(x) \cap N_1$, then $y_1, y_2\in N(x')$, and then $\{x, x', y_1, y_2\}$ induces either a diamond or a $K_4$ in $G$, a contradiction. So, $|N(x)\cap N_1| = 1$.

\smallskip
\no{$(iii)$}:  Since $C$ is a component of $[V(G)\sm (N_0\cup N_1)]$, by (i), there exists a vertex $y \in Y$ such that $N(y)\cap V(C) \neq \es$. We may assume that $y \in Y_1$. Now, we show that $N(y)\cap V(C) = V(C)$. Suppose not. Then there exist vertices $x, z \in V(C)$ such that $yx \in E$ and $yz \notin E$. Then since $C$ is connected, there exists a path joining $x$ and $z$ in $C$, say $P$. But, then $V(P) \cup \{y, v_2, v_3, v_4\}$ will induce a  $P_6$ in $G$, a contradiction. So (iii) holds. $\diamond$

\medskip
By (3)(iii), we see that for each $j \geq 3$, $N_j = \es$.  So, $V(G) = N_0 \cup W \cup Y \cup N_2$.

Also, by (3)(iii) and (R2),  $[N_2]$ is a union of $K_2$'s and $K_1$'s, and hence $[N_2]$ is bipartite. Let $(N_2', N_2'')$ be a bipartition of $[N_2]$ such that $N_2'$ is a maximal independent set of $N_2$.

\begin{enumerate}
\item[(4)] For each  $i \in [5]$, $i$ mod $5$, we have: (i) $[Y_i\cup Y_{i+2},  N(Y_{i+1}\cup Y_{i+3} \cup Y_{i+4}) \cap N_2'']= \es$, and (ii) $[Y_{i+1} \cup Y_{i+3},  N(Y_i \cup Y_{i+2}) \cap N_2''] = \es$.
\end{enumerate}

\no{\it Proof of $(4)$}. $(i)$: We prove for $i =1$. Suppose to the contrary that there exist adjacent vertices $y \in Y_1\cup Y_3$ and $x \in N(Y_{2}\cup Y_{4} \cup Y_{5}) \cap N_2''$. Since  $N_2'$ is a maximal independent set of $N_2$, there exists $x' \in N_2'$ such that $xx'\in E$. Also, since $x \in N(Y_{2}\cup Y_{4} \cup Y_{5}) \cap N_2''$, there exists a vertex $y' \in Y_{2}\cup Y_{4} \cup Y_{5}$ such that $xy' \in E$. But, then $\{y, y'\} \subseteq N(x) \cap N_1$, a contradiction to (3)(ii). So (i) holds.

\smallskip
\no{$(ii)$}: Similar to the proof of $(i)$. $\diamond$

\medskip
Now, by using the above properties, we prove the theorem in three cases as follows:
\medskip

\no{\bf Case 1}.  $W_i= \es$, for each $i\in [5]$.

Define $S_1 : = \{v_1, v_3\} \cup Y_1 \cup Y_3 \cup (N(Y_2\cup Y_4 \cup Y_5) \cap N_2'')$, $S_2 := \{v_2, v_4\} \cup Y_2 \cup Y_4 \cup (N(Y_1\cup Y_3) \cap N_2'')$,  $S_3:= \{v_5\} \cup Y_5$, and $S_4 := N_2'$. Then by (2)(iv) and (4),  $S_1, S_2$, $S_3$ and $S_4$ are independent sets. So, $(S_1, S_2, S_3, S_4)$ is a $4$-coloring of $G$.

\medskip
\no{\bf Case 2}. $W_i \neq \es$, for every $i \in [5]$.

By (2)(iii), $W_i$ is an independent set, for each $i$. Now, we define $S_1 := \{v_1, v_3\} \cup W_2 \cup Y_1 \cup Y_3 \cup (N(Y_2\cup Y_4 \cup Y_5) \cap N_2'')$, $S_2:= \{v_2, v_4\} \cup W_3 \cup Y_2 \cup Y_4 \cup (N(Y_1\cup Y_3) \cap N_2'')$, $S_4:= \{v_5\} \cup W_1 \cup Y_5$, and $S_5 := W_4 \cup W_5 \cup N_2'$. Then by the above properties, we see that $(S_1, S_2, S_3, S_4)$ is a $4$-coloring of $G$.

\medskip
\no{\bf Case 3}.  $W_i \neq \es$ and $W_{i-1}= \es$, for some $i \in [5]$, $i$ mod $5$.

Up to symmetry, we may assume that $i = 1$. Then by (2)(iii),  $W_3$ and $W_4$ are independent sets.

 (a) Suppose that $W_4= \es$. Then we define $S_1:= \{v_1, v_3\} \cup W_2' \cup Y_1 \cup Y_3 \cup (N(Y_2\cup Y_5) \cap N_2'')$, $S_2 := \{v_2, v_5\} \cup W_1' \cup Y_2 \cup Y_5 \cup (N(Y_1\cup Y_3 \cup Y_4) \cap N_2'')$, $S_3:= \{v_4\}\cup W_1''\cup W_2''\cup N_2'$, and $S_4 := W_3 \cup Y_4$. Then by the above properties,  $(S_1, S_2, S_3, S_4)$ is a $4$-coloring of $G$.

(b) Suppose $W_4 \neq \es$, then by (2)(iii), $W_1$ and $W_2$ are independent sets. Now, we define $S_1:= \{v_1, v_3\}\cup W_2 \cup Y_1\cup Y_3 \cup (N(Y_2\cup Y_5) \cap N_2'')$,  $S_2:= \{v_2, v_5\} \cup W_1 \cup Y_2 \cup Y_5 \cup (N(Y_1\cup Y_3 \cup Y_4) \cap N_2'')$, $S_3:= \{v_4\} \cup W_3 \cup Y_4$, and $S_4:= W_4 \cup N_2'$. Then by the above properties,  $(S_1, S_2, S_3, S_4)$ is a 4-coloring of $G$.

This completes the proof of the theorem.   \hfill{$\Box$}

\medskip
The bound given in Theorem~\ref{PDKB-4-col} is tight. For example, consider the duplicated graph $H_i$ obtained from $G_i$, $i \in \{2, 3\}$ (shown in Figure~\ref{ext}) by substituting each vertex indicated in circle by an independent set (of order $\geq 1$).   Then it is verified that both $H_2$ and $H_3$ are ($P_6$, diamond, $K_4$, bull)-free (using (R5)), and  $\chi(H_2) = \chi(H_3) = 4$. Note that the graph $G_2$ is a Greenwood-Gleason graph or a Clebsch graph, and the graph $G_3$ is the Mycielski 4-chromatic graph or a Gr\"otzsch graph.

 \begin{thm}\label{PDK-contains-Bull}
Let $G$ be a connected ($P_6$, diamond, $K_4$)-free graph that contains an induced bull. Then $G$ is $6$-colorable.
\end{thm}

\no{\it Proof.} Let $G$ be a connected ($P_6$, diamond, $K_4$)-free graph that contains an induced bull, say $H$ with vertex-set
$N_0 := \{v_1, v_2, v_3, v_4, v_5\}$, and  edge-set $\{v_1v_2,v_2v_3,v_3v_4,v_2v_5,v_3v_5\}$. For any integer $j\geq 1$, let $N_j$ denote the set $\{ y\in V(G) \mid d(y, N_0)=j \}$. Then we have the following:

\begin{enumerate}
\item[(1)] If $x\in N_1$, then $|N(x)\cap \{v_2, v_3, v_5\}|\leq 1$ and so $|N(x)\cap N_0|\leq 3$.

{\it Proof.} If  $|N(x)\cap \{v_2, v_3, v_5\}| \geq 2$, then $\{v_2, v_3, v_5, x\}$ induces a diamond or a $K_4$ in $G$, a contradiction. $\diamond$
\end{enumerate}

\noindent{By} (1) and since $G$ is diamond-free, we partition $N_1$ as follows. Let:
\begin{eqnarray*}
A_{i} &:=& \{x \in N_1 \mid  N(x)\cap N_0 = \{v_i\}\},  i \in \{1, \ldots, 5\}, \\
A_{jk} &:=& \{x \in N_1 \mid  N(x)\cap N_0 =\{v_i, v_{j}\}\}, j\in\{1, 4\} ~\mbox{and}~ k \in \{2, 3, 4, 5\}~ (j \neq k), \\
A_{p14} & := & \{x\in N_1 \mid  N(x)\cap N_0 = \{v_{p}, v_1, v_4\}\}, p \in \{2, 3, 5\}.
\end{eqnarray*}

\begin{enumerate}
\item[(2)]

By (R2) and (R3), for any $i \in \{1, \ldots, 5\}$, $j\in \{1, 4\}$, $k\in \{2, 3, 4, 5\}$   ($j \neq k$) and $p \in \{2, 3, 5\}$, we see that: $[A_i]$ is a union of $K_2$'s and $K_1$'s, and hence bipartite, and $[A_{jk}]$ and $[A_{p14}]$ are independent sets.
\end{enumerate}

\no{}For each $i\in \{1, \ldots, 5\}$, let $(A_i', A_i'')$ be a bipartition of $A_i$ such that $A_i'$ is a maximal independent set.

\begin{enumerate}
\item[(3)] We have either $A_1 = \es$ or $A_4 = \es$.

{\it Proof.} Suppose not. Let $x \in A_1$ and $y \in A_4$. Then since $\{x, v_1, v_2, v_3, v_4, y\}$ does not induce a $P_6$ in $G$, $xy \in E$. But then $\{v_1, x, y, v_4, v_3, v_5\}$ induces a $P_6$ in $G$, which is a contradiction. $\diamond$
\end{enumerate}

\begin{enumerate}
\item[(4)] We have the following:
\begin{enumerate}[label=(\roman*)]
\item $[A_1, A_{14} \cup A_{15}\cup A_{45}]$ is complete.
\item $[A_5, A_{15}\cup A_{45}]= \es$.
\item $[A_2'', A_{12} \cup A_{24} \cup A_{124}] = \es = [A_3'', A_{13}\cup A_{34}\cup A_{134}]$.
\item $A_2'' \cup A_{12} \cup A_{24} \cup A_{124}$ and $A_3'' \cup A_{13}\cup A_{34}\cup A_{134}$ are independent sets.
\item $[A_5'', A_{14} \cup A_{145}] = \es$.
\item $[A_{14} \cup A_{15}\cup A_{45}, A_{145}] = \es$.

\end{enumerate}
{\it Proof}. $(i)$: Suppose not. Then there exist vertices $x \in A_1$ and $y \in A_{14} \cup A_{15} \cup A_{45}$ such that $xy \notin E$. But, then
$\{x, v_1, y, v_3, v_4, v_5\}$ or $\{x, v_1, v_2, v_3, v_4, y\}$ induces a $P_6$ in $G$, a contradiction. So (i) holds.

$(ii)$: Suppose not. Then there exist vertices $x \in A_5$ and $y \in A_{15} \cup A_{45}$ such that $xy \in E$. But, then $\{y, x, v_1, v_2, v_3, v_4\}$ induces a $P_6$ in $G$, a contradiction. So, (ii) holds.

$(iii)$: Follows by the definitions of $A_2''$ and $A_3''$, and by (R1).

$(iv)$: Follows by (2), (iii), (R1), (R3), and (R4).

$(v)$:  $[A_5'', A_{145}] = \es$ follows by the definition of $A_5''$ and by(R1). Also, $[A_5'', A_{14}] = \es$. Otherwise, there exist
$x \in A_5''$ and $y \in A_{14}$ such that $xy \in E$. Since $A_5'$ is maximal, there exists $x' \in A_5'$ such that $xx' \in E$. But, then by (R1), $\{x', x, y, v_1, v_2, v_3\}$ induces a $P_6$ in $G$ which is a contradiction. So, $[A_5'', A_{14}] = \es$.

$(vi)$: Follows by (R3).
$\diamond$
\end{enumerate}

\no{}Define $S_1:= A_2'' \cup A_{12} \cup A_{24} \cup A_{124}$ and $S_2:= A_3'' \cup A_{13}\cup A_{34}\cup A_{134}$.

\begin{enumerate}
\item[(5)] If $C$ is a component of $[V(G)\sm(N_0\cup N_1)]$, then there exists a vertex $x \in N_1$ such that  $N(x)\cap V(C)= V(C)$.

{\it Proof}. Suppose not. Since $G$ is connected  there exists a vertex $x\in N_1$ such that $N(x)\cap V(C) \neq \es$ and $N(x) \cap V(C) \neq V(C)$. Then there exist vertices $y_1$ and $y_2$ in $C$ with $y_1y_2 \in E$ such that $y_1\in N(x)$ and $y_2\notin N(x)$. %By (1), $|N(x)\cap \lbrace v_2,v_3,v_5\rbrace |\leq 1$.
Now, if both $v_1$ and $v_4$ are neighbors of $x$ or if both $v_1$ and $v_4$ are non-neighbors of $x$, then by (1), $N_0 \cup \{x, y_1, y_2\}$ induces a $P_6$ in $G$, a contradiction.
So, we may assume, up to symmetry that $xv_1 \in E$ and $xv_4 \notin E$. Now, if $xv_2 \in E$, then by (1), $\{y_2, y_1, x, v_2, v_3, v_4\}$ induces a $P_6$ in $G$, a contradiction. So, $xv_2 \notin E$. Then since $\{y_2, y_1, x, v_1, v_2, v_5\}$ does not induce a $P_6$ in $G$, $xv_5 \in E$. But, then $\{y_2, y_1, x, v_5, v_3, v_4\}$ induces a $P_6$ in $G$, a contradiction. So, (5) holds. $\diamond$
\end{enumerate}

\no{}By (5) and by (R2), $[N_2]$ is a union of $K_2$'s and $K_1$'s, and hence bipartite. Let $(N_2', N_2'')$ be a bipartition of $N_2$ such that $N_2'$ is a maximal independent set of $N_2$.

Also, since $G$ is connected, by (5) and by the definition of $N_j$'s, $N_j= \es$, for all $j\geq 3$. Thus, $V(G) = N_0 \cup N_1 \cup N_2$.

\begin{enumerate}
\item[(6)] $[A_1\cup A_4\cup A_{15}\cup A_{45}, N_2]= \es$.

{\it Proof}. Suppose not. Then there exist vertices $x \in A_1 \cup A_4 \cup A_{15}\cup A_{45}$ and $y \in N_2$ such that $xy \in E$. But, then $\{x, y, v_1, v_2, v_3, v_4\}$ induces a $P_6$ in $G$, a contradiction.
$\diamond$
\end{enumerate}

\no{}By (4:(iv)), it enough to show that $G - (S_1 \cup S_2)$ is $4$-colorable, and we do this in two cases using (3).

\medskip

\no{\bf Case 1}. {\it Suppose that $A_1 \cup A_4 \neq \es$.}

By (3) and by symmetry we may assume that $A_1 \neq \es$ and $A_4 = \es$. Then:
\begin{clm}\label{A1-neq-es} The following hold:

\begin{enumerate}[label=(\roman*)]
\item  $[A_1, A_5] = \es$.

%\item $[A_1, A_{15}],[A_1,A_{45}]$ are complete.

\item $A_{14} \cup A_{15} \cup A_{45}$ is an independent set.

\item $[A_3, N_2]= \es = [A_5, N_2]$.

\end{enumerate}
\end{clm}

\no{\it Proof of Claim~$\ref{A1-neq-es}$}.
$(i)$: Suppose not. Then there exist vertices $x\in A_1$ and $y \in A_5$ such that $xy \in E$. But, then $\{y, x, v_1, v_2, v_3, v_4\}$ induces a $P_6$ in $G$, a contradiction. So (i) holds.  $\diamond$

\no{$(ii)$}: Suppose not. Then there exist adjacent vertices, say $x$ and $y$ in $A_{14} \cup A_{15} \cup A_{45}$. Since $A_1 \neq \es$, let $z \in A_1$. Then by (4:(i)), $[\{z\}, A_{14} \cup A_{15} \cup A_{45}]$ is complete. But, then $\{z, x, v_4, y\}$ or $\{z, x, v_5, y\}$  induce a diamond in $G$ or
$\{z, x, v_1, y\}$ induces a $K_4$ in $G$, a contradiction. So (ii) holds.
$\diamond$

\no{$(iii)$}:
Suppose not. Since $A_1 \neq \es$, let $z \in A_1$.

Let $x \in A_3$ and $y \in N_2$ be such that $xy \in E$.  Then since $\{z, v_1, v_2, v_3, x, y\}$ does not induce a $P_6$ in $G$, $xz \in E$. But, then $\{v_5, v_2, v_1, z, x, y\}$ induces a $P_6$ in $G$, a contradiction. So, $[A_3, N_2]= \es$.

Again, let $x\in N_2$ and $y\in A_5$ be such that $xy \in E$.  By $(6)$ and $(i)$, $xz, yz\notin E$. But, then $\{z, v_1, v_2, v_5, y, x \}$ induces a $P_6$ in $G$, a contradiction. So, $[A_5, N_2]= \es$. $\diamond$

Hence we have proved Claim~\ref{A1-neq-es}.  $\Diamond$

\medskip

\no{Now,} let us define $S_3:= A_2' \cup \{v_1, v_5\}$;  $S_4:= A_3' \cup N_2' \cup \{v_2, v_4\}$; $S_5:= A_1' \cup A_5' \cup N_2'' \cup \{v_3\}$; and $S_6:= A_1'' \cup A_5'' \cup A_{14}\cup A_{15} \cup A_{45} \cup A_{145}$. Then by the above claim and by the above properties, we see that $(S_3, S_4, S_5, S_6)$ is a $4$-coloring of $G -(S_1 \cup S_2)$.

\medskip
\no{\bf Case 2}. {\it Suppose that $A_1 \cup A_4 = \es$.}

\medskip
If $[A_2 \cup A_3 \cup A_5, N_2]\neq \es$, then we find a suitable bull with $A_1\neq \es$, and we proceed as in Case~1 to get a $4$-coloring of $G-(S_1 \cup S_2)$.  Also, if $[A_{14}, A_{15}]$ is not complete, then there exist $u\in A_{14}$  and $v\in A_{15}$ such that $uv \notin E$. Now,  $\{v_2, v_3, v_5, v_4, v\}$ induces a bull with $A_1 \neq \emptyset$. So, we proceed as in Case~1 to get a $4$-coloring of $G-(S_1 \cup S_2)$. By symmetry, the same holds if $[A_{14}, A_{45}]$ is not complete.  So, we may assume that $[A_2 \cup A_3 \cup A_5, N_2]= \es$, and that $[A_{14}, A_{15}]$  and $[A_{14},A_{45}]$ are complete.

  We define $S_3:= N_2'\cup A_2' \cup \{v_1, v_5\}$ and $S_4:= N_2''\cup A_3' \cup \{v_2, v_4\}$. Clearly $S_3$ and $S_4$ are independent sets. Now it is enough to show that $G\setminus (S_1\cup S_2\cup S_3\cup S_4)$ is bipartite.
 Then:
\begin{clm} \label{empty set in N1}
$A_{14} =\emptyset$ or $A_{15}=\emptyset$ or $A_{45}=\emptyset$.
\end{clm}
\no{\it Proof of Claim~$\ref{empty set in N1}$.}  Suppose not, let $x\in A_{14}, y\in A_{15}$ and $z\in A_{45}$. Then since $\{y, v_1, v_2, v_3, v_4, z\}$ does not induce a $P_6$ in $G$, $yz \in E$. But then since $[A_{14}, A_{15}]$ and $[A_{14},A_{45}]$ are complete, $\{x, y, z, v_5\}$ induces a diamond in $G$, a contradiction. So the claim holds.  $\diamond$

\medskip
Now, we define
\begin{equation*}
S_5 :=
\begin{cases}
 A_5' \cup A_{15}\cup \{v_3\}, \  \mbox{if}~ A_{14}=\emptyset \ \mbox{or}~  A_{45}=\emptyset, \\
 A_5' \cup A_{45} \cup \{v_3\}, \  \mbox{if}~ A_{15}=\emptyset,
\end{cases}
\end{equation*}
and
\begin{equation*}
S_6 :=
\begin{cases}
A_5'' \cup A_{14} \cup A_{45} \cup A_{145}, ~\mbox{if}~ A_{14}=\emptyset ~ \mbox{or} \ A_{45}=\emptyset. \\
A_5'' \cup A_{14} \cup A_{145}, ~ \mbox{if}~ A_{15}=\emptyset.
\end{cases}
\end{equation*}
 Then by the above claim and by (4), we see that $(S_5, S_6)$ is a bipartition of $G\setminus (S_1 \cup S_2\cup S_3\cup S_4)$.

This completes the proof of the theorem. \hfill$\Box$

\medskip
The bound given in Theorem~\ref{PDK-contains-Bull} is tight. For example, consider the graph $G$ which is isomorphic to the complement of the 16-regular \emph{Schl\"afli graph} on 27 vertices. As mentioned earlier,  $G$ is ($P_6$, diamond)-free,  $\chi(G) = 6$ and $\omega(G) = 3$. Also, it is verified that $G$ contains a bull.

\begin{thm}
Let $G$ be a ($P_6$,  diamond, $K_4$)-free graph. Then $G$ is $6$-colorable.
\end{thm}
\no{\it Proof}. Follows by Theorems~\ref{PDKB-4-col} and \ref{PDK-contains-Bull}. \hfill{$\Box$}

\end{document}